\documentclass[amstex,12pt, amssymb]{article}

\usepackage{mathtext}
\usepackage[cp1251]{inputenc}
\usepackage[T2A]{fontenc}
\usepackage[dvips]{graphicx}
\usepackage{amsmath}
\usepackage{amssymb}
\usepackage{amsxtra}
\usepackage{latexsym}
\usepackage{ifthen}

\newcounter{lemma}[section]

\newcounter{corollary}[section]

\newcounter{remark}[section]

\newcounter{theorem}[section]

\newcounter{proposition}[section]

\numberwithin{equation}{section}

\begin{document}

\markboth{\centerline{VLADIMIR RYAZANOV}} {\centerline{INFINITE
DIMENSION OF SOLUTIONS FOR DIRICHLET PROBLEM II}}

\author{{VLADIMIR RYAZANOV}}

\title{{\bf INFINITE DIMENSION OF SOLUTIONS \\ FOR DIRICHLET PROBLEM II}}

\maketitle

\large \begin{abstract} It is proved that the space of solutions of
the Dirichlet problem for the harmonic functions in the unit disk
with nontangential boundary limits $0$ a.e. has the infinite
dimension.
\end{abstract}


\bigskip

By the well--known Lindel\"of maximum principle, see e.g. Lemma 1.1
in \cite{GM}, it follows the uniqueness theorem for the Dirichlet
problem in the class of bounded harmonic functions $u$ on the unit
disk $\mathbb D = \{ z\in\mathbb{C}: |z|<1\}$. In general there is
no uniqueness theorem in the Dirichlet problem for the Laplace
equation even under zero boundary data. In comparison with the
previous arXiv versions and \cite{RI}, here we give more elementary
examples and constructions of solutions.

\medskip

Many such nontrivial solutions $u$ for the Laplace equation can be
given by the {\bf Poisson-Stiltjes integral}
$$
u(z)\ =\ \frac{1}{2\pi}\ \int\limits_{0}\limits^{2\pi}P_r(\vartheta
-t)\ d\Phi(t)\ ,\ \ \ z=re^{i\vartheta}, \ r<1\
,\quad\quad\quad\quad\quad\quad (1)
$$
with an arbitrary {\bf singular function} $\Phi:[0,2\pi]\to\mathbb
R$, i.e., where $\Phi$ is of bounded variation and $\Phi^{\prime}=0$
a.e., where we use the standard notation for the {\bf Poisson
kernel}
$$
P_r(\Theta)\ =\ \frac{1-r^2}{1-2r\cos\Theta +r^2}\ , \ r<1\ .
$$

Indeed, $u$ in (1) is harmonic for every function
$\Phi:[0,2\pi]\to\mathbb R$ of bounded variation and by the Fatou
theorem, see e.g. Theorem I.D.3.1 in \cite{Ko},
$u(z)\to\Phi^{\prime}(\Theta)$ as $z\to e^{i\Theta}$ along any
nontangential path whenever $\Phi^{\prime}(\Theta)$ exists. Thus,
$u(z)\to0$ as $z\to e^{i\Theta}$ for  a.e. $\Theta\in[0,2\pi]$ along
any nontangential paths for every singular function $\Phi$.

\medskip

{\bf Example 1.} The simplest example of such kind is given by
nondecreasing step-like data $\Phi_{\vartheta_0}$ with values $0$
and $2\pi$ and with the jump at $\vartheta_0\in(0,2\pi)$:
$$
u(z)\ =\ P_r(\vartheta -\vartheta_0)\ =
\frac{1-r^2}{1-2r\cos(\vartheta -\vartheta_0) +r^2}\ ,\ \ \
z=re^{i\vartheta}, \ r<1\ . \quad (2)
$$
We directly see that $u(z)\to0$ as $z\to e^{i\Theta}$ for all
$\Theta\in(0,2\pi)$ except $\Theta=\vartheta_0$.

Note that the function $u$ is harmonic in the unit disk $\mathbb D$
because
$$
u(z)\ =\ {\rm {Re}}\ \frac{\zeta_0+z}{\zeta_0-z}\ =\
\frac{1-|z|^2}{1-2\, {\rm {Re}}\, z\overline{\zeta_0} +|z|^2}\
,\quad\quad \zeta_0=e^{i\vartheta_0}\ ,\ z\in\mathbb D\ , \quad\quad
(3)
$$
where the function $w=g(z)=g_{\zeta_0}(z)\colon =
(\zeta_0+z)/(\zeta_0-z)$ is analytic (conformal) in $\mathbb D$ and
maps $\mathbb D$ onto half-plane ${\rm {Re}}\, w>0$, $g(0)=1$,
$g(\zeta_0)=\infty$.

\medskip

{\bf Example 2.} The second natural example is given by the formula
(1) with $\Phi(t)=\varphi(t/2\pi)$ where $\varphi:[0,1]\to[0,1]$ is
the well--known {\bf Cantor function}, see e.g. \cite{DMRV} and
further references therein.

\medskip

The formula (2) gives a continual set of such examples. Furthermore,
one can prove the following result.

\medskip

{\bf Theorem 1.} {\it The space of all harmonic functions in
$\mathbb D$ with nontangential limit $0$ at every point of
$\partial\mathbb D$ except a countable collection of points in
$\partial\mathbb D$ has the infinite dimension.}

\medskip

\begin{proof} Indeed, let us consider the sequence of
functions of the form (3)\ :
$$
u_n(z)\ =\ {\rm {Re}}\ \frac{\zeta_n+z}{\zeta_n-z}\ =\
\frac{1-|z|^2}{1-2\, {\rm {Re}}\, z\overline{\zeta_n} +|z|^2}\
,\quad\quad \zeta_n=e^{i\vartheta_n}\ ,\ z\in\mathbb D\ ,
$$
where $$\vartheta_n=\pi(2^{-1}+\ldots +2^{-n}),\quad n=1,2,\ldots
$$
and denote by ${\cal H}_1$ the class of all series $u=\sum\gamma_n
u_n$ whose sequences of coefficients $\gamma=\{\gamma_n\}$ belong to
the space $l^1$ with the norm
$\|\gamma\|=\sum\limits_{n=1}\limits^{\infty}|\gamma_n|<\infty$.
Note that ${\cal H}_1$ consists of harmonic functions, see, e.g.,
Theorem I.3.1 in \cite{Go}, because
$$
0\ <\ u_n(z)\ <\ \frac{1+|z|}{1-|z|}\,\quad \forall\ n=1,2,\ldots\
,\ z\in\mathbb D\ .
$$

Note also that each function $u\in {\cal H}_1$ has nontangential
limit $0$ at every point $\zeta\in\partial\mathbb D$ except the
points $\zeta_0=-1=e^{i\vartheta_0}$, $\vartheta_0=\pi$, and
$\zeta_n$, $n=1,2,\ldots$. Indeed, let $\zeta = e^{i\Theta}$,
$\Theta\in(0,2\pi)$, $\zeta\ne\zeta_n$, $n=0,1,2,\ldots$. Then,
applying the formula (2), we have the estimate
$$
u_n(z)\ \leq\ \frac{1-r^2}{4r\sin^2\frac{\Theta -\vartheta_n}{4}}\
\leq\ C(1-r)\ ,\ \ \ z=re^{i\vartheta},
$$
for all points $z$ belonging to a sector $|\vartheta -
\Theta|<c(1-r)$ and for all $r$ which are close enough to $1$ where
$C< \infty$ does not depend on $n=1,2,\ldots$. Thus,
$$
|u(z)|\ \le\ C\, \|\gamma\|\, (1-r)\ \to\ 0\quad\quad {\rm as}\quad
r\to 1\ ,\ \ \ z=re^{i\vartheta},
$$
in any sector $|\vartheta - \Theta|<c(1-r)$.

Now, let us show that $u_n$, $n=1,2,\ldots $, form a basis in the
space ${\cal H}_1$ with the locally uniform convergence in $\mathbb
D$ which is metrizable.

Indeed, firstly, $u=\sum\limits_{n=1}\limits^{\infty}\gamma_nu_n\ne
0$ if $\gamma\ne 0$. Really, let us assume that $\gamma_n\ne 0$ for
some $n=1,2,\ldots $. Then $u\ne 0$ because $u(z)\to\infty$ as
$z=re^{i\vartheta_n}\to
 e^{i\vartheta_n}$. The latter follows because
$$u_n(re^{i\vartheta_n})\ =\ \frac{1+r}{1-r}\ \to\ \infty \quad\quad {\rm as}\quad
r\to 1\ ,$$ and by the previous item
$$
|{\tilde {u}}(re^{i\vartheta_n})|\ \le\ C\, \|{\gamma}\|\, (1-r)\
\to\ 0\quad\quad {\rm as}\quad r\to 1\ ,
$$
where ${\tilde {u}}=u-\gamma_nu_n$.

\bigskip

Secondly, $u^*_m=\sum\limits_{n=1}\limits^{m}\gamma_nu_n\to u$
locally uniformly in $\mathbb D$ as $m\to\infty$. Indeed, elementary
calculations give the following estimate of the remainder term
$$
|u(z)-u^*_m(z)|\ \leq\ \frac{1+r}{1-r}\ \cdot
\sum\limits_{n=m+1}\limits^{\infty}|\gamma_n|\ \to\ 0\ \ \ \ \ \ \ \
\ {\mbox{as}}\ \ \ \ \ m\to\infty\quad (4)
$$
in every disk $\mathbb D(r)=\{ z\in\mathbb C: |z|\leq r\}$, $r<1$.
\end{proof}$\Box$

\bigskip

{\bf Corollary 1.} {\it Given a measurable function $\varphi
:\partial\mathbb D\to \mathbb R$, the space of all harmonic
functions $u:\mathbb D\to\mathbb R$ with the limits
$\lim\limits_{z\to\zeta}u(z) = \varphi(\zeta)$ for a.e.
$\zeta\in\partial\mathbb D$ along nontangential paths has the
infinite dimension.}

\medskip

Indeed, the existence at least one such a harmonic function $u$
follows from the known Gehring theorem in \cite{Ge}. Combining this
fact with Theorem 1, we obtain the conclusion of Corollary 1.

\bigskip

{\bf Remark 1.} In view of Lemma 3.1 in \cite{ER}, one can similarly
prove the more refined result on harmonic functions than in
Corollary 1  with respect to logarithmic capacity instead of the
measure of the length on $\partial\mathbb D$.

\medskip

Moreover, the statements on the infinite dimension of the space of
solutions can be extended to the Riemann-Hilbert problem because the
latter is reduced in the papers \cite{ER} and \cite{RI} to the
corresponding two Dirichlet problems.

\medskip

Note also that harmonic functions $u$ found in Theorem 1 themselves
cannot be represented in the form of the Poisson integral with any
integrable function $\Phi:[0,2\pi]\to\mathbb R$ because this
integral would have nontangential limits $\Phi$ a.e., see e.g.
Corollary IX.9.1 in \cite{Go}. Consequently, $u$ do not belong to
the classes $h_p$ for any $p>1$, see e.g. Theorem IX.2.3 in
\cite{Go}.

However, the functions $u\in {\cal H}_1$ in the proof of Theorem 1
have the representation as the Poisson-Stiltjes integral (1) with
$\Phi =\sum\gamma_n \Phi_{\vartheta_n}$ where
$\Phi_{\vartheta_n}:[0,2\pi]\to\mathbb R$ are nondecreasing
step-like functions with values $0$ and $2\pi$ with jumps at the
points $\vartheta_n$, $n=1,2,\ldots$. Thus, $\Phi$ is of bounded
variation and hence ${\cal H}_1\subset h_1$, see e.g. Theorem IX.2.2
in \cite{Go}.

\medskip

{\bf Problem 1.} It remains the open question whether the basis of
the space of all such singular solutions of the Dirichlet problem
for the Laplace equation has the power of the continuum.

\medskip
\noindent
{\bf Vladimir Illich Ryazanov,}\\
Institute of Applied Mathematics and Mechanics,\\
National Academy of Sciences of Ukraine,\\
74 Roze Luxemburg Str., Donetsk, 83114, Ukraine,\\
vl.ryazanov1@gmail.com

\end{document}